\input amstex
\documentstyle{amsppt}
\magnification=\magstep1 \NoRunningHeads
\input pictexwd.tex

\title
INFINITE MEASURE PRESERVING TRANSFORMATIONS
WITH RADON MSJ
\endtitle

\author
Alexandre I. Danilenko
\endauthor

\email
alexandre.danilenko@gmail.com
\endemail

\address
 Institute for Low Temperature Physics
\& Engineering of National Academy of Sciences of Ukraine, 47 Lenin Ave.,
 Kharkov, 61164, UKRAINE
\endaddress
\email alexandre.danilenko\@gmail.com
\endemail

\abstract
We introduce  concepts of  Radon MSJ and  Radon disjointness for infinite Radon measure preserving homeomorphisms of the locally compact Cantor space.
We construct an uncountable family of pairwise  Radon disjoint infinite Chacon like transformations.
Every such transformation is Radon strictly ergodic, totally ergodic, asymmetric (not isomorphic to its inverse), has  Radon MSJ and possesses  Radon  joinings whose ergodic components are not joinings.
\endabstract

 \loadbold

\NoBlackBoxes

\endtopmatter

\document

\head 0. Introduction
\endhead

This paper (inspired by a recent progress in \cite{JaRoRu}) is about property of {\it minimal self-joinings} (MSJ) for {\it infinite measure preserving} transformations.
The concept of MSJ for probability preserving maps was introduced in \cite{Rud}
as a powerful tool to construct systems with certain prescribed dynamical properties
(see also \cite{dJRu}, \cite{Ru}, \cite{Da2} and references therein for further developments).
Some subsequent works are devoted to  extension of MSJ to infinite measure preserving and nonsingular systems: \cite{RudSi}, \cite{SiWi}, \cite{dJSi}, \cite{JaRoRu}, etc.
However such extensions are only partial due to principal obstacles  such as the following ones:
\roster
\item"---"
given a continuous transformation of a Polish space with a recurrent point, there exist uncountably many ergodic pairwise orthogonal nonsingular (as well as infinite invariant) measures \cite{Gli}, \cite{Ef},
 \item"---"
ergodic components of a non-ergodic joining  need not be joinings  \cite{Aa1},
 \item"---"
given an ergodic  probability preserving transformation $S$ and $\lambda\in(0,1)$, there is en ergodic 2-fold self-joining of $S$ of Krieger's type $III_\lambda$ \cite{RudSi}.
\endroster
To bypass these obstacles the authors of the aforementioned papers have to select special classes of measures invariant or quasi-invariant under Cartesian powers of the transformations under question.
For instance, only {\it boundedly finite} infinite measures were considered  in \cite{JaRoRu} and mostly {\it rational} nonsingular joinings were studied in \cite{RudSi}.
We also go this way.

 We first develop elements of the general  theory of joinings for  ergodic measure preserving homeomorphisms $T$ of a locally compact non-compact Cantor space $X$ endowed with an infinite  Radon measure $\mu$.
This is not a restriction because  every ergodic infinite (non-atomic) $\sigma$-finite measure preserving transformation on a standard Borel space is measure theoretically isomorphic to a Radon strictly ergodic homeomorphism on $X$ \cite{Yu}.
Considering only Radon joinings and self-joinings of ergodic Radon dynamical systems,  we introduce  notions of {\it Radon MSJ} and {\it Radon disjointness}.
Using the topological nature of $X$, $T$ and $\mu$ and Hopf ratio ergodic theorem we can define {\it generic points} for the system.

In the remaining part of the paper we present examples of transformations with Radon MSJ.
First we consider  the {\it infinite Chacon transformation} $T$.
It was introduced in 1997 by Adams, Friedman and Silva
 (see \cite{AdFrSi, Section 2}) as
an infinite measure preserving counterpart of the
 classical Chacon map.
In a recent paper \cite{JaRoRu} Janvresse, Roy and de la Rue showed that
this transformation has MSJ in the class of so-called  boundedly finite (in fact, Radon) measures\footnote{This extends the well known result that the classical Chacon has MSJ \cite{dJRaSw}.} and that
there exist  ergodic  infinite $T^{\times d}$-invariant  measures whose marginals are singular with respect to the original $T$-invariant measure.
Such measures are called {\it weird} in \cite{JaRoRu}.
We  provide  different (shorter and more algebraic) proofs of these results in Theorem~5.2.
While in \cite{AdFrSi} and \cite{JaRoRu}, $T$ is a Borel bijection (with countably many discontinuities) of a co-countable subset of the infinite interval $[0,\infty)$ furnished with Lebesgue measure,
in this paper $T$ appears as a Radon strictly ergodic homeomorphism of  $X$.
In addition to what was done in \cite{JaRoRu} we
\roster
\item"$\bullet$" describe explicitly the recurrent points for $T\times T$ (Proposition~4.3),
\item"$\bullet$" describe the weird $T^{\times d}$-invariant measures (for each $d>1$) as
{\it quasi-graph} measures (see Section 3.2),
\item"$\bullet$" introduce an uncountable family of {\it infinite Chacon like} Radon strictly ergodic transformations $T^{\omega}$, $\omega\in\{0,1\}^\Bbb N$ (in Section~6),
\item"$\bullet$" show that each $T^\omega$ is totally ergodic,  has  Radon MSJ and is not isomorphic to the inverse $(T^\omega)^{-1}$ (Theorem~6.1, Corollary~6.2),
\item"$\bullet$" prove that $T^\omega$ and $T^{\omega'}$ are isomorphic if the infinite sequences $\omega$ and $\omega'$ are tail equivalent
and that $T^\omega$ and $T^{\omega'}$ are  Radon disjoint otherwise (Theorem~6.1),
\item"$\bullet$" give an example of nonergodic Radon 2-fold joining of $T$ such that almost all ergodic components of this joining are  conservative quasi-graphs whose coordinate projections are singular to the original $T$-invariant measure (Example~5.4).
\endroster

The outline of the paper is as follows.
In Section~1 we introduce some basic notions for Radon dynamical systems such as recurrent points and  generic points and discuss ergodic decompositions.
Radon joinings and Radon disjointness for general Radon dynamical systems
are studied in Section~2.
In Section~3 we recall the $(C,F)$-construction of rank-one actions for Abelian groups.
We introduce their  quasi-graph invariant measures and study their properties.
The $(C,F)$-construction is used in Section~4 to define the infinite Chacon transformation $T$.
The recurrent points for the Cartesian square of $T$ are  explicitly described there.
In Section~5 we show that $T$ has  Radon MSJ.
As a  byproduct we describe all $T^{\times d}$-invariant measures on $X^d$.
An example of a Radon joining of $T$ whose ergodic components are not joinings is also given in this section.
In Section~6 we introduce and study the infinite Chacon like transformations.
In the final Section~7 we discuss possible generalizations of the class of dynamical systems to which the main results of this paper extend.
Some open problems are also stated there.

\subsubhead{Acknowledgements} \endsubsubhead
I thank the anonymous referee for useful remarks.

\head 1. Basic definitions: recurrent~points,~generic~points, ergodic decomposition
\endhead

Let $X$ be a locally compact non-compact Cantor space and let $T$ be an aperiodic homeomorphism of $X$.

\definition{Definition 1.1}
A point $x\in X$ is called {\it $T$-recurrent} if for each neighborhood $U$ of $x$ there is an integer $n\ne 0$ such that $T^nx\in U$.
\enddefinition

We are interested in $T$-invariant (non-negative) $\sigma$-finite Borel measures on $X$.
Since $X$ is non-compact,  $T$-invariant probability measures on $X$ need not exist.
On the other hand,
it follows from \cite{Gli} and \cite{Sc} that if there exists a $T$-recurrent point then there are uncountably many non-equivalent  infinite $T$-invariant $\sigma$-finite  Borel measures on $X$.
In this work our main concern is Radon measures.
We recall that a Borel measure $\lambda$ on $X$ is called {\it Radon} if it is finite on every compact subset of $X$.
Of course, every Radon measure is $\sigma$-finite.
The converse is not true.
Denote by $\Cal M_{\text{Ra}}(X)$ the cone of all Radon measures on $X$.
It is a Polish space in the $*$-weak topology.
Let $\Cal M_{\text{Ra}}(X,T)$ stand for the  subset of $T$-invariant Radon measures on $X$ and let $\Cal M^e_{\text{Ra}}(X,T)$ denote the subset of ergodic $T$-invariant Radon measures on $X$.
The two subsets are Borel.
Given a pair of measures $\lambda,\lambda'\in \Cal M^e_{\text{Ra}}(X,T)$,
we have that either $\lambda=c\lambda'$ for some real $c>0$ or $\lambda
\perp\lambda'$.
If $\Cal M_{\text{Ra}}(X,T)$ consists of a single ray we say that $T$ is {\it Radon uniquely ergodic}.
If, in addition, $T$ is minimal then we say that $T$ is {\it Radon strictly ergodic}.

\definition{Definition 1.2}
Let  $\lambda\in \Cal M^e_{\text{Ra}}(X,T)$.
A point $x\in X$ is called {\it generic} for $\lambda$ if for each pair of compact open subsets $K_1,K_2\subset X$ such that $\lambda(K_2)>0$, the following limit exists and is given by
$$
\lim_{n\to\infty}\frac{\sum_{a_n\le j\le b_n}1_{K_1}(T^jx)}{\sum_{a_n\le j\le b_n}1_{K_2}(T^jx)}= \frac{\lambda(K_1)}{\lambda(K_2)}
$$
whenever $a_n\le 0\le b_n$ and $b_n-a_n\to+\infty$.
\enddefinition

If $\lambda$ is non-atomic (i.e. the dynamical system $(X,\lambda,T)$ is conservative) then every generic point is recurrent.
Denote by $\Cal G(\lambda)$ the set of all generic points for $\lambda$.
Then  $\Cal G(\lambda)$ is a Borel subset of $X$ and $T\Cal G(\lambda)=\Cal G(\lambda)$ (pointwise, not only almost everywhere).
By the Hopf ratio ergodic theorem \cite{Aa},  $\lambda$-almost every point of $X$ is generic for $\lambda$\footnote{This theorem was stated there only for conservative systems. It is also true for dissipative (ergodic, Radon) ones with a trivial proof.}.
If $\lambda'$ is another ergodic  $T$-invariant Radon  measure then $\Cal G(\lambda)=\Cal G(\lambda')$  if ${\lambda'}=c\lambda$  for some $c>0$ and $\Cal G(\lambda)\cap\Cal G(\lambda')=\emptyset$ otherwise.

Let  $\eta\in\Cal M_{\text{Ra}}(X,T)$.
 Then there is a probability measure $\kappa$ on $\Cal M^e_{\text{Ra}}(X,T)$ such that
$\eta=\int_{\Cal M^e_{\text{Ra}}(X,T)}\lambda\,d\kappa(\lambda)$.
This integral is called an {\it ergodic decomposition} of $\eta$.
It is not unique.

\head 2. Radon self-joinings and Radon disjointness
\endhead

Let $X$ be a locally compact Cantor space, $G$ a discrete countable infinite Abelian group,
$T=(T_g)_{g\in G}$ a continuous action of $G$ on $X$ and $\mu$ a $T$-invariant Radon measure on $X$.
We call the triple $(X,T,\mu)$ {\it Radon dynamical system}.
From now on we assume that $\mu$ is ergodic.

\definition{Definition 2.1} Let $d>1$.
A Radon $((T_g)^{\times d})_{g\in G}$-invariant measure on $X^d$ whose coordinate projections (marginals) are all equivalent\footnote{We note that these projections are not necessarily $\sigma$-finite.}  to $\mu$ is called a {\it d-fold Radon self-joining of $T$}.
The set of all   Radon $d$-fold self-joinings of $T$ will be denoted by $J_{d,\text{Ra}}(T)$.
The subset of  ergodic    Radon $d$-fold self-joinings of $T$ will be denoted by $J_{d,\text{Ra}}^{e}(T)$.
\enddefinition

For instance, the product $\mu\otimes\cdots\otimes\mu(d \text{ times})$ belongs to $J_{d,\text{Ra}}(T)$.

Denote by $C(T)$ the centralizer of $T$, i.e. the group of all invertible
$\mu$-nonsingular transformations commuting with $T_g$ for each $g\in G$.
If $S_1,\dots,S_{d-1}\in C(T)$ then there are  constants $c_i>0$ such that $\mu\circ S_i=c_i\mu$ for all $i=1,\dots,d$\footnote{Because $T$ is ergodic.}.
The measure $\mu_{S_1,\dots,S_{d-1}}$ on $X^d$, given by
$$
\mu_{S_1,\dots,S_{d-1}}(A_1\times\cdots\times A_d ):=\mu(A_1\cap S_1^{-1}A_2\cap\cdots\cap S_{d-1}^{-1}A_d)
$$
 for all Borel subsets $A_1,\dots,A_d\subset X$, belongs to $J_{d,\text{Ra}}^{e}(T)$.
Since $\lambda_{S_1,\dots,S_{d-1}}$ is supported on the graph of the map $X\ni x\mapsto (S_1x,\dots, S_{d-1}x)\in X^{d-1}$, it is called a {\it graph-joining} of $T$.

\definition{Definition 2.2}
If for each $\lambda\in J_{d,\text{Ra}}^{e}(T)$ there is a partition of $\{1,\dots,d\}$ into subsets $J_1,\dots,J_k$
such that $\lambda$ splits into direct product of its marginals $\lambda_i$ on $X^{J_i}$, each $\lambda_i$ is (up to multiplicative constant) a graph-joining of $T$ if $\# J_i>1$ or $\lambda_i=\mu$ if $\# J_i=1$, $i=1,\dots,k$,  and $C(T)=\{T_g\mid g\in G\}$ then
we say that  $T$ has the property of {\it Radon  MSJ$_d$}.
If $T$ has Radon MSJ$_d$ for each $d>1$ then we say that $T$ has {\it Radon MSJ}.
\enddefinition

We note that if $\lambda$ is a non-ergodic Radon self-joining of $T$ then the ergodic components of $\lambda$ need not be self-joinings of $T$.
For instance, there exist  Radon strictly ergodic homeomorphism $R$ of $X$ with a $T$-invariant Radon measure $\mu$ such that the system $(X\times X,\mu\times\mu, R\times R)$ is totally dissipative \cite{Da1, Theorem~0.1(4)}, i.e. almost every ergodic component of the Radon 2-fold self-joining $\mu\times\mu$ of $R$ is an $(R\times R)$-orbit equipped with the ``counting'' purely atomic measure.
Purely atomic measures are not self-joinings of $R$ because their coordinate projections are also atomic and hence singular to $\mu$ which is continuous.

We also need the following analogue of the Furstenberg concept of disjointness for probability preserving systems \cite{Fu}.

\definition{Definition 2.3} Let $(X,(T_g)_{g\in G},\mu)$ and $(Y, (S_g)_{g\in G},\nu)$ be two ergodic Radon dynamical systems.
We say that they are {\it Radon disjoint} if $\mu\times\nu$ is the only (up to a multiplicative constant) ergodic Radon $(T_g\times S_g)_{g\in G}$-invariant measure on $X\times Y$ whose coordinate projections are equivalent to
$\mu$ and $\nu$ respectively.
\enddefinition

Of course, if $T$ and $S$ are Radon disjoint then they are not isomorphic.

\head 3. $(C,F)$-construction of rank-one Abelian actions
\endhead
\subhead 3.1. $(C,F)$-nomenclature
\endsubhead
For a detailed exposition of the $(C,F)$-construction we refer the reader to \cite{Da1} and \cite{Da3}.

Let $G$ be a discrete countable infinite Abelian group.
Let $(F_n)_{n\ge 0}$  and $(C_n)_{n\ge 1}$ be two sequences of finite subsets in $G$ such that for each $n>0$,
\roster
\item"(I)"
$F_0=\{0\}$,  $\# C_n>1$,
\item"(II)"
$F_0\subset F_1\subset \cdots$ with $\bigcup_{n\ge 0}F_n=G$,
\item"(III)"
$F_n+F_n+C_{n+1}\subset F_{n+1}$,
\item"(IV)"
$(F_n+c)\cap (F_n+c')=\emptyset$ for all $c,c'\in C_{n+1}$ with $c\ne c'$.
\endroster
We let $X_n:=F_n\times C_{n+1}\times C_{n+2}\times\cdots$ and endow this set of infinite sequences with the infinite product topology.
Then $X_n$ is a compact Cantor (i.e. totally disconnected  perfect metric) space.
The  mapping
$$
X_n\ni (f_n,c_{n+1},c_{n+2},\dots)\mapsto(f_n+c_{n+1}, c_{n+2},\dots)\in X_{n+1}
$$
is a topological embedding
of $X_n$ into $X_{n+1}$.
From now on we can consider $X_n$ an a subset of $X_{n+1}$.
Let $X$ denote  the inductive limit of the sequence $(X_n)_{n\ge 0}$ furnished with these embeddings.
  Then $X$ is a well defined locally compact Cantor  space.
We call it the {\it $(C,F)$-space associated with the sequence} $(C_n,F_{n-1})_{n\ge 1}$.
Given a subset $A\subset F_n$, we let
$$
[A]_n:=\{x=(f_n,c_{n+1},\dots)\in X_n\mid f_n\in A\}
$$
and call this set an  {\it $n$-cylinder} in $X$.
It is open and compact in $X$.
It is easy to verify that
$$
\align
[A]_n\cap[B]_n &=[A\cap B]_n, \quad
[A]_n\cup[B]_n =[A\cup B]_n\quad \text{and}\\
 [A]_n &=[A+C_{n+1}]_{n+1}
 \endalign
 $$
  for all $A,B\subset F_n$ and  $n\ge 0$.
For brevity, we will write $[f]_n$ for $[\{f\}]_n$, $f\in F_n$.
We note that the collection of all cylinders coincides with the family of  all compact open subset  in $X$.

Let $\Cal R$ denote the {\it tail equivalence relation} on $X$.
This means that the restriction of $\Cal R$ to $X_n$ is the tail equivalence relation on $X_n$ for each $n\ge 0$.
We note that $\Cal R$ is {\it minimal}, i.e. the $\Cal R$-class of every point  is dense in $X$.
There exists a {\it unique} $\sigma$-finite $\Cal R$-invariant Borel measure $\mu$ on $X$ such that
$\mu(X_0)=1$.
It is  Radon.
Moreover,  $\mu$ is strictly positive on every non-empty open subset.
We note that the {\it $\Cal R$-invariance}\footnote{$\mu$ is called $\Cal R$-invariant if $\mu$ is invariant under each Borel transformation whose graph is contained in $\Cal R$.} of $\mu$ is equivalent to the following property:
$$
\mu([f]_n)=\mu([f']_n)\quad\text{ for all }f,f'\in F_n, n\ge 0.
$$
This and the normalization condition on $\mu$ imply that
$$
\mu([A]_n)=\frac{\# A}{\# C_1\cdots\# C_n} \quad\text{for each subset }A\subset F_n,\  n>0.
$$
We call $\mu$ the {\it $(C,F)$-measure associated with $(C_n, F_{n-1})_{n\ge 1}$}.
It is infinite if and only if\footnote{In view of (I)--(IV), the sequence $(\frac{\# F_n}{\# C_1\cdots\#C_n})_{n=1}^\infty$ is non-decreasing and bounded by $1$ from below.}
$$
\lim_{n\to\infty}\frac{\# F_n}{\# C_1\cdots\#C_n}=\infty.
\tag3-1
$$
It is easy to see that  $\mu$  is {\it $\Cal R$-ergodic}, i.e. each Borel $\Cal R$-saturated subset of $X$ is either $\mu$-null or $\mu$-conull.
We now define an action of $G$ on $X$.
Given $g\in G$, let
$$
X_n^g:=\{(f_n,c_{n+1},c_{n+2}\dots)\in X_n\mid g+f_n\in F_n\}.
$$
Then $X_n^g$ is a compact open subset of $X_n$ and $X_n^g\subset X_{n+1}^g$.
In view of (III) and~(II),  $X_n\subset X_{n+1}^g$ eventually in $n$.
Hence  $\bigcup_{n\ge 0}X_n^g=X$.
Therefore, given $x\in X$ and $g\in G$, there is $n>0$ such that if we consider an expansion $x=(f_n,c_{n+1},\dots)$ of $x$ in $X_n$ then $g +f_n\in F_n$.
We now let
$$
T_gx:=(g+f_n, c_{n+1},\dots)\in X_n\subset X.
$$
It is standard to verify that
\roster
\item"(i)" the map $T_g:X\ni x\mapsto T_gx\in X$ is a homeomorphism of $X$ and
\item"(ii)"$T_gT_{g'}=T_{g+g'}$ for all $g,g'\in G$.
\endroster
Thus $T:=(T_g)_{g\in G}$ is a continuous action of $G$ on $X$.
We call it {\it the $(C,F)$-action of $G$ associated with the sequence $(C_n,F_{n-1})_{n\ge 1}$}.
It is free.
  The $T$-orbit equivalence relation is  $\Cal R$.
Indeed, if $(x,x')\in\Cal R$ then there is $n>0$ such that $x,x'\in X_n$.
Consider now the expansions $x=(f_n,c_{n+1},\dots)$ and $x'=(f'_n,c_{n+1}',\dots)$ of $x$ and $x'$ in $X_n$.
Then $x'=T_gx$, where $g:=f_n'-f_n+\sum_{j> n}(c_j'-c_j)$.
It follows that  $T$ preserves  $\mu$.
We call $(X,\mu,T)$ the {\it dynamical system associated with} $(C_n,F_{n-1})_{n\ge 1}$.
It is Radon strictly ergodic\footnote{We also note that $T$ has rank one along the sequence $(F_n)_{n=1}^\infty$ (in the measure theoretical sense, see \cite{Da1}).}.

We will need the following standard lemma (see \cite{Da4, Lemma~1.1(ii)} for the proof of a more general fact).

\proclaim{Lemma 3.1} Let $g\in G$.
If there is $\delta>0$ such that for each $n\in\Bbb N$ and $v,w\in F_n$, there is a subset $A\subset[v]_n$ and $m\in\Bbb Z$ such that $T_{mg}A\subset[w]_n$
and $\mu(A)\ge \delta\mu([v]_n)$  then the transformation $T_g$ is ergodic.
\endproclaim

\subhead 3.2. Quasi-graph  invariant Radon measures
\endsubhead
In this subsection we introduce an important class of Radon $((T_g)^{\times d})_{g\in G}$-invariant measures on $X^d$ which are ``close'' to graph-joinings but whose marginals can be singular to $\mu$.
Such measures can appear as ergodic components of non-ergodic Radon $d$-fold self-joinings of $T$ (see \S 5 and \S 6).
Throughout this subsection we will assume that the following condition holds:
\roster
\item"(V)"
$F_n+C_{n+1} +F_n-F_n\subset F_{n+1}$  for each $n>0$.
\endroster
This condition is stronger than (III).
Let $J$ a subset of $\Bbb N$.
For each $n\not\in J$, fix an element   $\alpha_n\in C_n$.
Thus we obtain  a sequence   $\alpha:=(\alpha_n)_{n\not\in J}$.
We now define a new $(C,F)$-sequence $(C_n^\alpha,F_{n-1}^\alpha)_{n\ge 1}$ by setting
$F^\alpha_n:=F_n$ and
$$
C_n^\alpha:=\cases
C_n, &\text{if }n\in J\\
\{\alpha_n\},&\text{if }n\not\in J.
\endcases
$$
Then  $(C_n^\alpha,F_{n-1}^\alpha)_{n\ge 1}$ satisfies (I)--(V) except for the fact that the cardinality of $C_n^\alpha$ can be  $1$. 
This happens if and only if $n\not\in J$.
Nevertheless the associated $(C,F)$-action $T^\alpha=(T^\alpha_g)_{g\in G}$ on a locally compact $(C,F)$-space $X^\alpha$ equipped with an invariant  $\sigma$-finite  measure $\mu^\alpha$ is well defined (via the same inductive construction process described in \S 3.1).
The space $X^\alpha$ is  totally disconnected.
It is either perfect (i.e. Cantor)  or purely discrete and countable.
The latter happens if and only if $J$ is finite.
We note that $\mu^\alpha$ is finite if and only if $\mu$ is finite and $J$ is cofinite\footnote{I.e. the complement $\Bbb N\setminus J$ is finite.}.
Given $A\subset F_n$, the corresponding cylinder in $X^\alpha$ will de denoted by $[A]_n^\alpha$  (in order to not confuse it with the cylinder $[A]_n$ in $X$).
We now define an embedding $\phi_\alpha$ of $X^\alpha$ into $X$ by setting for each $n>0$ and a point $x=(f_n,c_{n+1},c_{n+2},\dots)\in X^\alpha_n\subset X^\alpha$,
$$
\phi_\alpha(x):=(f_n,c_{n+1},c_{n+2},\dots)\in X_n\subset X.
$$
It is straightforward to verify that $\phi_\alpha$ is a well-defined, one-to-one continuous equivariant\footnote{I.e. $\phi^\alpha \circ T_g^\alpha=T_g\circ \phi^\alpha$ for each $g\in G$.} map.
Hence the image $\phi^\alpha(X^\alpha)$ of $X^\alpha$ is a $T$-invariant dense $F_\sigma$-subset of $X$.
Now given $m>0$ and a subset $A\subset F_m$, take a point $x\in[A]_m\cap\phi_\alpha(X^\alpha)$.
If we write the expansion $x=(f_m,c_{m+1},c_{m+2},\dots)$ of $x$  in $X_m$ then $f_m\in A$ and $c_l=\alpha_l$ for all $l\not\in J$  that are greater than some $n>m$.
It follows that
$$
[A]_m\cap \phi_\alpha(X^\alpha)=\bigsqcup_{g\in \sum_{J\not\ni l>m}(C_l-\alpha_l)} T_g\phi_\alpha([A]_m^\alpha).\tag3-2
$$
We have the following trichotomy: 
\roster
\item"$\bullet$"
if $J$ is cofinite then $\mu_\alpha\circ \phi_\alpha^{-1}$  is a multiple of $\mu$ and hence $\phi_\alpha$ is an isomorphism (both topological and measure theoretical) of $(X^\alpha,\mu^\alpha,T^\alpha)$ onto  $(X,\mu,T)$;
\item"$\bullet$"
if $J$ is finite then $\mu_\alpha\circ \phi_\alpha^{-1}\perp\mu$ and $(X^\alpha,\mu^\alpha,T^\alpha)$ is an ergodic totally dissipative dynamical system;
\item"$\bullet$"
if $J$ is neither finite nor cofinite then $\mu_\alpha\circ \phi_\alpha^{-1}\perp\mu$ and $(X^\alpha,\mu^\alpha,T^\alpha)$ is an ergodic conservative dynamical system.
\endroster

Consider now another sequence $\beta=(\beta_n)_{n\not\in J}$ with $\beta_n\in C_n$ for each $n\not\in J$.
Let $(X^\beta,\mu^\beta,T^\beta)$ be the dynamical system associated with $(C_n^\beta,F_{n-1}^\beta)_{n\ge 1}$, where $C_n^\beta$ and $F_{n-1}^\beta$ are defined in a similar way as $C_n^\alpha$ and $F_{n-1}^\alpha$.
Then there is a natural measure preserving topological equivariant isomorphism $\phi_{\alpha,\beta}$ of $(X^\alpha,\mu^\alpha,T^\alpha)$  onto $(X^\beta,\mu^\beta,T^\beta)$.
To define $\phi_{\alpha,\beta}$ we take a point
 $x\in X^\alpha_{n-1}$ for some $n>0$ and write the expansion of $x$ in $X^\alpha_n$:
$x=(f_n,c_{n+1},c_{n+2},\dots)\in F_n\times C_{n+1}^\alpha\times C_{n+2}^\alpha\cdots$.
Then $f_n\in F_{n-1}+C_n^\alpha$.
We now set 
$$
\phi_{\beta,\alpha}(x):=(\widetilde f_n,\widetilde c_{n+1},\widetilde c_{n+2},\dots),
$$
where $\widetilde f_n:=f_n+\sum_{J\not\ni j\le n}(\beta_j-\alpha_j)$ and  for each $l>n$,
$$
\widetilde c_l:=
\cases
c_l,&\text{if }l\in J\\
\beta_l,&\text{if }l\not\in J.
\endcases
$$
We note that  if $n\in J$ then $\sum_{J\not\ni j\le n}(\beta_j-\alpha_j)\in F_{n-1}-F_{n-1}$
and hence $\widetilde f_n\in f_n+F_{n-1}-F_{n-1}\subset F_{n-1}+C_n^\alpha+F_{n-1}-F_{n-1}$.
If $n\not\in J$ then $\sum_{J\not\ni j\le n}(\beta_j-\alpha_j)\in \beta_n-\alpha_n+F_{n-1}-F_{n-1}$ and hence $\widetilde f_n\in F_{n-1}+C_n^\alpha+\beta_n-\alpha_n+F_{n-1}-F_{n-1}=F_{n-1}+C_n^\beta+F_{n-1}-F_{n-1}$.
Therefore in every of the two possible cases
$\widetilde f_n\in F_{n}$ in view of (V).
Hence $\phi_{\alpha,\beta}(x)\in X_n^\beta$.
Thus $\phi_{\alpha,\beta}$ is well defined.
Moreover, for each subset $A\subset F_{n-1}+C_n^\alpha$, 
$$
\phi_{\beta,\alpha}([A]_n^\alpha)=[\widetilde A]^\beta_n,\text{ where }
\widetilde A:=A+\sum_{J\not\ni j\le n}(\beta_j-\alpha_j).
\tag3-3
$$
It is straightforward  to verify that $\phi_{\beta,\alpha}$ is as claimed.
It is easy to see that the converse $\phi_{\beta,\alpha}^{-1}$ to $\phi_{\beta,\alpha}$ is $\phi_{\alpha,\beta}$.
We now define a map $\psi_{\alpha,\beta}$ from  $X^\alpha$ to $X\times X$ by setting
$$
\psi_{\alpha,\beta}(x):=(\phi_\alpha(x),\phi_\beta(\phi_{\beta,\alpha}(x))).
$$
It is one-to-one, continuous and  $\psi_{\alpha,\beta}\circ T^\alpha_g =(T_g\times T_g)\circ\psi_{\alpha,\beta}$ for all $g\in G$.

\proclaim{Theorem 3.2} 
Suppose that for each $l\not\in J$,
$$
(C_l-\alpha_l+F_{l-1}^\bullet)\cap(C_l-\beta_l+F_{l-1}^\bullet)=F_{l-1}^\bullet,\tag3-4
$$
where $F_{l-1}^\bullet:=F_{l-2}+C_{l-1}-F_{l-2}-C_{l-1}$.
Then the following are satisfied.
\roster
\item"\rom(i)"
The subset
  $\psi_{\alpha,\beta}(X^\alpha)$ is closed in $X\times X$.
Hence $\psi_{\alpha,\beta}(X^\alpha)$ is locally compact in the induced topology.
\item"\rom(ii)" 
The topological dynamical system  $(\psi_{\alpha,\beta}(X^\alpha), (T_g\times T_g)_{g\in G})$ is Radon strictly ergodic.
\item"\rom(iii)"
$\mu^\alpha\circ\psi^{-1}_{\alpha,\beta}$ is the only (up to scaling) $T_1\times T_1$-invariant Radon measure supported on $\psi_{\alpha,\beta}(X^\alpha)$.
\item"\rom(iv)" 
Given $a,b\in F_m$, represent the cylinder $[a]_m\times[b]_m$ in $X\times X$ as the infinite product
$\{(a,b)\}\times (C_{m+1}\times C_{m+1})\times(C_{m+2}\times C_{m+2})\times\cdots$.
Then 
$$
{\mu^\alpha\circ\psi^{-1}_{\alpha,\beta}\restriction([a]_m\times [b]_m)}
={\mu^\alpha([a]^\alpha_m)}\cdot\widehat\delta_{a',b'}\otimes\bigotimes_{n>m}\gamma_n,
$$
where $\gamma_m$ is a probability on $C_m\times C_m$ given by
$$
\gamma_n:=
\cases
\frac 1{\# C_n}\sum_{c\in C_n}\delta_c\times\delta_c,&\text{if $n\in J$}\\
\delta_{(\alpha_n,\beta_n)}, &\text{if $n\not\in J$},
\endcases
$$
$a'=a+\sum_{J\not\ni j\le m}\alpha_j$ and $b'=b+\sum_{J\not\ni j\le m}\beta_j$.\footnote{Here and below $\delta_c$ or $\delta_{(a,b)}$ denotes the probability measure supported at a single point $\{c\}$  or $\{(a,b)\}$ respectively while $\widehat\delta_{a',b'}$ denotes the Kronecker delta.}
\endroster
\endproclaim

\demo{Proof}
(i) If $J$ is cofinite then it is straightforward to verify that $\psi_{\alpha,\beta}(X_\alpha)$ is
 the the graph of  $T_h$, where $h:=\sum_{j\not\in J}(\beta_j-\alpha_j)$.
 The graph is closed in $X\times X$.
 
 Consider now another particular case where $J$ is finite.
 Then it is straightforward to verify that $\psi_{\alpha,\beta}(X_\alpha)$ is the $(T_g\times T_g)_{g\in G}$-orbit of a point $(x,y)\in X\times X$ such that $x=(f_m,\alpha_{m+1},\alpha_{m+2},\dots)\in X_m$
 and $y=(r_m,\beta_{m+1},\beta_{m+2},\dots)\in X_m$ for some $m\ge 0$ and elements $f_m,r_m\in F_m$.
 It suffices to show that given $l>m$ and subsets $A,B\subset F_l$, the intersection of the orbit of $(x,y)$ with the cylinder $[A]_l\times[B]_l$ is finite.
 Indeed, $(T_gx,T_gy)\in[A]_l\times[B]_l$ if and only if
 $$
 \align
 &g\in A-f_m-\alpha_{m+1}-\cdots-\alpha_l+\sum_{i>l}(C_i-\alpha_i)\quad\text{and}\\
 & g\in B-r_m-\beta_{m+1}-\cdots-\beta_l+\sum_{i>l}(C_i-\beta_i).
  \endalign
 $$
 Hence $g\in  F_{l+1}^\bullet+\sum_{i>l+1}(C_i-\alpha_i)$ and $g\in  F_{l+1}^\bullet+\sum_{i>l+1}(C_i-\beta_i)$.
 It follows from~\thetag{3-4} that $( F_{l+1}^\bullet+\sum_{i>l+1}(C_i-\alpha_i))\cap( F_{l+1}^\bullet+\sum_{i>l+1}(C_i-\beta_i))=F_{l+1}^\bullet$.
 It remains to note that $F_{l+1}^\bullet$ is finite.

Suppose  now that $J$ is neither finite nor cofinite.
Fix $n>0$ and subsets $A\subset F_{n-1}+C_n$ and $B\subset F_{n-1}+C_n^\beta$.
Applying~\thetag{3-2}  we obtain that
$$
\align
\psi^{-1}_{\alpha,\beta}&([A]_n\times[B]_n)\subset \phi_\alpha^{-1}([A]_n)\cap\phi_{
\beta,\alpha}^{-1}(\phi_\beta^{-1}([B]_n))\\
&=
\left(\bigsqcup_{g\in \sum_{J\not\ni l>n}(C_l-\alpha_l)} T_g^\alpha[A]_n^\alpha\right)\cap
\phi_{
\alpha,\beta}\left(\bigsqcup_{g\in \sum_{J\not\ni l>n}(C_l-\beta_l)} T_g^\beta[ B]_n^\beta\right)
\\
&=\left(\bigsqcup_{J\not\ni l>n}[A+C_{n+1}+\cdots+C_l]_l^\alpha\right)\cap
\phi_{
\alpha,\beta}\left(\bigsqcup_{J\not\ni l>n}[ B+C_{n+1}+\cdots+C_l]_l^\beta\right)
\\
&=\bigsqcup_{J\not\ni l>n}([A+C_{n+1}+\cdots+C_l]_l^\alpha\cap
\phi_{
\alpha,\beta}([ B+C_{n+1}+\cdots+C_l]_l^\beta)).
\endalign
$$
Given $l\not\in J$, let $l'$ be the least integer such that $l'>l$ and $l'\not\in J$.
Then $C_k^\alpha=C_k^\beta=C_k$ if $l<k<l'$, $C_{l'}^\alpha=\{\alpha_{l'}\}$,  $C_{l'}^\beta=\{\beta_{l'}\}$ and hence
$$
\align
[A+C_{n+1}+\cdots+ C_l]^\alpha_l&=[A'+\alpha_{l'}]^\alpha_{l'}\quad\text{and}\\
[B+C_{n+1}+\cdots+ C_l]^\beta_l&=[B'+\beta_{l'}]^\beta_{l'},
\endalign
$$
where $A':=A+C_{n+1}+\cdots+ C_{l'-1}\subset F_{l'-1}$ and 
$B':=B+C_{n+1}+\cdots+ C_{l'-1}\subset F_{l'-1}$.
Utilizing this and \thetag{3-3} we obtain that\footnote{In the formulas below, $l'=l'(l)$.}
$$
\align
\psi^{-1}_{\alpha,\beta}([A]_n\times[B]_n)
&
\subset 
\bigsqcup_{J\not\ni l>n}\left( [A'+\alpha_{l'}]_{l'}^\alpha\cap
\bigg[ B'+\beta_{l'}+ \sum_{J\not\ni j\le l'}(\alpha_j-\beta_j)           \bigg]_{l'}^\alpha\right)\\
&=\bigsqcup_{J\not\ni l>n}
\left( \bigg[\bigg(A'\cap
\bigg(B'+ \sum_{J\not\ni j\le l}(\alpha_j-\beta_j)\bigg)\bigg)    +\alpha_{l'}       \bigg]_{l'}^\alpha\right).
\endalign
$$
It follows from \thetag{3-4} that for all subsets $D,D'\subset F_n^\bullet$,
$$
\bigg(D+\sum_{j=n+1}^{l'-1}C_{j}-\sum_{n<j<l'}^{j\not\in J}\alpha_j\bigg)\cap
\bigg(D'+\sum_{j=n+1}^{l'-1}C_{j}-\sum_{n<j<l'}^{j\not\in J}\beta_j\bigg)=(D\cap D')+\sum_{n<j<l'}^{j\in J}C_j.
$$
Therefore we have
$$
\aligned
A'\cap
&\bigg(B'+ \sum_{J\not\ni j\le l}(\alpha_j-\beta_j)\bigg)  \\
&=
\bigg(\bigg(A+\sum_{j=n+1}^{l'-1} C_j-\sum_{j<l'}^{j\not\in J}\alpha_j\bigg)\cap\bigg(B +\sum_{j=n+1}^{l'-1} C_j-\sum_{j<l'}^{j\not\in J}\beta_j\bigg)\bigg)
+\sum_{j<l'}^{j\not\in J}\alpha_j\\
&=
\bigg(\bigg(\widehat A+\sum_{j=n+1}^{l'-1} C_j-\sum_{n<j<l'}^{j\not\in J}\alpha_j\bigg)\cap\bigg(\widehat B +\sum_{j=n+1}^{l'-1} C_j-\sum_{n<j<l'}^{j\not\in J}\beta_j\bigg)\bigg)
+\sum_{j<l'}^{j\not\in J}\alpha_j\\
&=(\widehat A\cap\widehat B)+ \sum_{n<j<l'}^{j\in J}C_j
+\sum_{j\le n}^{j\not\in J}\alpha_j+\sum_{n<j<l'}^{j\not\in J}\alpha_j\\
&=(A\cap\widetilde B)+ \sum_{n<j<l'}C_j^\alpha,
\endaligned
$$
where 
$\widehat A:=A-\sum_{J\not\in j\le n}\alpha_j\subset F_{n-1}+C_n-F_{n-1}-C_n=F_n^\bullet$,
$\widehat B:=B-\sum_{J\not\in j\le n}\beta_j\subset F_n^\bullet$
and
$\widetilde B:=B +\sum_{J\not\ni j\le n}(\alpha_j-\beta_j)$.
Hence 
$$
\psi^{-1}_{\alpha,\beta}([A]_n\times[B]_n) \subset\bigsqcup_{J\not\ni l>n}\bigg[(A\cap \widetilde B) +\sum_{n<j\le l'} C_j^\alpha\bigg]^\alpha_{l'}= [A\cap\widetilde B]_n^\alpha.
$$
On the other hand, it is straightforward to verify that the converse inclusion holds.
Therefore, 
$$
\psi^{-1}_{\alpha,\beta}([A]_n\times[B]_n) =[A\cap \widetilde B]^\alpha_{n}.\tag3-5
$$
It follows that the $\psi_{\alpha,\beta}$-preimage of $[A]_n\times[B]_n$   is compact  in $X^\alpha$.
Hence $[A]_n\times[B]_n$ does not intersect the set $\overline{\psi_{\alpha,\beta}(X_\alpha)}\setminus
\psi_{\alpha,\beta}(X_\alpha)$.
Since $J$ is not finite then we take
 $A:=B:=F_{n-1}+C_n$  for $n\in J$ and obtain that the union
$$
\bigcup_{n \in J}([F_{n-1}+C_n]_{n}\times[F_{n-1}+C_n]_n)=\bigcup_{n\in J} (X_{n-1}\times X_{n-1})=X\times X
$$
 does not intersect  $\overline{\psi_{\alpha,\beta}(X_\alpha)}\setminus
\psi_{\alpha,\beta}(X_\alpha)$.
Hence  $\psi_{\alpha,\beta}(X_\alpha)$
 is closed.
 \comment
 If $J$ is finite then we take $A:=F_{n-1}+\beta_n$ and $B:=F_{n-1}+\beta_n$ 
 for $n\not\in J$ and obtain that the set $[F_{n-1}+\beta_n]_n\times[F_{n-1}+\beta_n]_n$ does not intersect 
 $\overline{\psi_{\alpha,\beta}(X_\alpha)}\setminus
\psi_{\alpha,\beta}(X_\alpha)$.
Since the latter set is invariant under $(T_g\times T_g)_{g\in  G}$, it follows that the set
$$
\bigcup_{g\in C_n-\beta_n}(T_g[F_{n-1}+\beta_n]_n\times T_g[F_{n-1}+\beta_n]_n)=[C]
$$
\endcomment

(ii)
It follows from (i) that $\psi$ is an equivariant  topological isomorphism of the dynamical system $(X^\alpha, T^\alpha)$ onto $(\psi_{\alpha,\beta}(X^\alpha), (T_g\times T_g)_{g\in G})$.
Since the former system is Radon strictly ergodic, so is the latter.  

(iii) is obvious now. 

(iv) 
Let $m,a,b$ be the same as in the statement of~(iv).
Suppose that $J$ is infinite.
Take $n>m$ with $n\in J$ and elements $c_j,\widehat c_j\in C_j$ whenever $m<j\le n$.
Then by~\thetag{3-5},
$$
\multline
\psi_{\alpha,\beta}^{-1}
\bigg(
\bigg[a+  \sum_{m<j\le n}c_j\bigg]_n
\times\bigg[b+\sum_{m<j\le n}\widehat c_j\bigg]_n
\bigg)
\\
=\bigg[\bigg\{a+\sum_{m<j\le n}c_j\bigg\}\cap
\bigg\{b+\sum_{m<j\le n}\widehat c_j+\sum_{J\not\ni j\le n}(\alpha_j-\beta_j)\bigg\}\bigg]_n^\alpha.
\endmultline
$$
Therefore, by \thetag{3-4}, the intersection in the above formula is non-empty if and only if $a=b+ \sum_{J\not\ni j\le m}(\alpha_j-\beta_j)$,  $c_j=\widehat c_j$ if $J\ni j$
and $c_j=\alpha_j$ and $\widehat c_j=\beta_j$ if $J\not\ni j$.
Hence
$$
\mu^\alpha\bigg(\psi^{-1}_{\alpha,\beta}
\bigg(
\bigg[a+  \sum_{j=m+1}^nc_j\bigg]_n
\times\bigg[b+\sum_{j=m+1}^n\widehat c_j\bigg]_n
\bigg)\bigg)=
\mu^\alpha([a]_m^\alpha)\cdot \widehat\delta_{a',b'}\prod_{j=m+1}^n\gamma_j(c_j,\widehat c_j),
$$
where $a',b'$ and $\gamma_j$ were introduced in the statement of (iv).
\qed
\enddemo

\definition{Definition 3.3} Given a subset $J\subset \Bbb N$, two sequences $\alpha=(\alpha_n)_{n\not\in J}$ and $\beta=(\beta_n)_{n\not\in J}$ with $\{\alpha_n,\beta_n\}\subset C_n$ for each $n\not\in J$ and an element $h\in G$,
we call the measure $(\mu^\alpha\circ\psi^{-1}_{\alpha,\beta})\circ ( I\times T_h^{-1})$ the {\it quasi-graph} measure.
\enddefinition

It follows from Theorem~3.2 that if \thetag{3-4} is satisfied then $(\mu^\alpha\circ\psi^{-1}_{\alpha,\beta})\circ ( I\times T_h^{-1})\in\Cal M^e_{\text{Ra}}(X\times X, (T_g\times T_g)_{g\in G})$.
The projection of this measure to the first coordinate is equivalent to $\mu^\alpha$ and
the projection of this measure to the second coordinate is equivalent to $\mu^\beta$.
The conditional measures corresponding to the disintegrations of 
$(\mu^\alpha\circ\psi^{-1}_{\alpha,\beta})\circ ( I\times T_h^{-1})$
over $\mu^\alpha$ (or over $\mu^\beta$) are delta-measures almost everywhere.
 The  quasi-graph Radon measure is a 2-fold self-joining of $T$ if and only if  $\mu^\alpha\sim\mu^\beta\sim\mu$.
 This happens if and only if  $J$ is cofinite. 
 In this case the quasi-graph measure is a graph-joining of $T$\footnote{The quasi-graph invariant Radon measures which are not graph-joinings are analogues of the so-called {\it weird measures} considered in \cite{JaRoRu} in the case when $T$ is an infinite Chacon transformation.}.

In a similar way we may define quasi-graph  invariant Radon measures for higher Cartesian powers of $T$. 
Indeed, given $d>1$
and sequences $\alpha^{(i)}=(\alpha^{(i)})_{n\not\in J}$, $i=1,\dots,d$, such that $\alpha_n^{(i)}\in C_n$ for all $i=1,\dots,d$, we define a map $\psi_{\alpha^{(1)},\dots,\alpha^{(d)}}$ from $X^{\alpha^{(1)}}$  to the $d$-th Cartesian power of $X$ by setting
$$
\psi_{\alpha^{(1)},\dots,\alpha^{(d)}}(x):=(\phi_{\alpha^{(1)}}(x),\phi_{\alpha^{(2)}}(\phi_{\alpha^{(2)},\alpha^{(1)}}(x)),\dots,
\phi_{\alpha^{(d)}}(\phi_{\alpha^{(d)},\alpha^{(1)}}(x))).
$$
This map is one-to-one, continuous and equivariant, i.e. $\psi_{\alpha^{(1)},\dots,\alpha^{(d)}}\circ (T^{\alpha^{(1)}})_1=(T_1)^{\times d}\circ \psi_{\alpha^{(1)},\dots,\alpha^{(d)}}$.
The following theorem is a $d$-fold analogue of Theorem~3.2.
It can be proved in a similar way.
Therefore we leave the proof  to the reader.

\proclaim{Theorem 3.4} Suppose that for each $l\not\in J$, there are  $u_l,v_l\in\{1,\dots,d\}$ such that
$$
(C_l-\alpha_l^{(u_l)}+ F_{l-1}^\bullet)\cap (C_l-\alpha_l^{(v_l)}+F_{l-1}^\bullet)=F_{l-1}^\bullet,
\tag3-6
$$
where $F_{l-1}^\bullet$ is same as in the statement of Theorem~3.2.
Then the following are satisfied.
\roster
\item"\rom(i)"
The subset $\psi_{\alpha^{(1)},\dots,\alpha^{(d)}}(X^{\alpha^{(1)}})$ is closed in $X^d$.
Hence it is locally compact in the induced topology.
\item"\rom(ii)" 
The topological dynamical system  $(\psi_{\alpha^{(1)},\dots,\alpha^{(d)}}(X^{\alpha^{(1)}}), ((T_g)^{\times d})_{g\in G})$ is Radon strictly ergodic.
\item"\rom(iii)"
$\mu^{\alpha^{(1)}}\circ\psi_{\alpha^{(1)},\dots,\alpha^{(d)}}^{-1}$ is the only (up to scaling) $((T_g)^{\times d})_{g\in G}$-invariant Radon measure supported on $\psi_{\alpha^{(1)},\dots,\alpha^{(d)}}(X^{\alpha^{(1)}})$.
\item"\rom(iv)" 
Given $a^{(1)},\dots,a^{(d)}\in F_m$, represent the cylinder $[a^{(1)}]_m\times\cdots\times[a^{(d)}]_m$ in $X^d$ as the infinite product
$\{(a^{(1)},\dots,a^{(d)})\}\times C_{m+1}^d\times C_{m+2}^d\times\cdots$.
Then 
$$
\multline
\mu^{\alpha^{(1)}}\circ\psi_{\alpha^{(1)},\dots,\alpha^{(d)}}^{-1}\restriction([a^{(1)}]_m\times\cdots\times[a^{(d)}]_m)\\
={\mu^{\alpha^{(1)}}([a]^{\alpha^{(1)}}_m)}\prod_{i=1}^d\widehat\delta_{a^{(1)}_\bullet,a^{(i)}_\bullet}\otimes\bigotimes_{n>m}\gamma_n,
\endmultline
$$
where $\gamma_n$ is a probability on $C_n\times C_n$ given by
$$
\gamma_n:=
\cases
\frac 1{\# C_n}\sum_{c\in C_n}\bigotimes_{i=1}^d\delta_c,&\text{if $n\in J$}\\
\delta_{(c_n^{(1)},\dots,c^{(d)}_n)}, &\text{if $n\not\in J$},
\endcases
$$
$a^{(i)}_\bullet:=a^{(i)}+\sum_{J\not\ni j\le m}\alpha_j^{(i)}$ and elements $c_n^{(i)}\in C_n$, $i=1,\dots,d$, are such that $c_n^{(u_n)}=\alpha_n^{(u_n)}$ and $c_n^{(v_n)}=\alpha_n^{(v_n)}$
for each $n>m$.
\endroster
\endproclaim

By analogy with Definition~3.3 we now give the following definition.

\definition{Definition 3.5} Given $\alpha^{(1)},\dots,\alpha^{(d)}$ as above and let $h_1,\dots,h_{d-1}$ be elements of $G$.
The measure $\mu^{\alpha^{(1)}}\circ\psi_{\alpha^{(1)},\dots,\alpha^{(d)}}^{-1}\circ(I\times T_{h_1}^{-1}\times\cdots\times T_{h_{d-1}}^{-1})$ is called a   {\it quasi-graph} measure.
\enddefinition

It follows from Theorem~3.4 that if \thetag{3-6} is satisfied for $\alpha^{(1)},\dots,\alpha^{(d)}$ then the corresponding quasi-graph measure is Radon, $((T_g)^{\times d})_{g\in G}$-invariant and  ergodic.
In this case a quasi-graph Radon measure is purely atomic if and only if $J$ is finite.
On the other hand, a quasi-graph Radon measure is a $d$-fold Radon joining ($d$-fold graph-joining indeed) if and only if $J$ is cofinite.

\head 4. Infinite Chacon transformation~and~recurrent~points for its Cartesian square
\endhead

To define the infinite Chacon transformation we will use the $(C,F)$-construction.
From now on $G=\Bbb Z$.
We first define recurrently a  sequence of positive integers  $(h_n)_{n=0}^\infty$ by setting
$h_0:=1$ and $h_{n+1}:=6h_n+1$ for $n>0$\footnote{In \cite{AdFrSi} and \cite{JaRoRu}, it was assumed that $h_{n+1}=6h_n+2$. To achieve the ``symmetric'' version of the infinite Chacon   we use one spacer less at each step of the inductive construction. See Section~7 for a discussion about other possibilities to add extra spacers.}.
For each $n\ge 0$, we now let
$$
\align
F_n&:=\left\{-\frac{h_n-1}2,-\frac{h_n-3}2,\dots,\frac{h_n-1}2\right\}\qquad\text{ and }\\
C_{n+1}&:=\{-h_n,0,h_n+1\}.
\endalign
$$
The sequence $(C_n,F_{n-1})_{n\ge 1}$ satisfies (I)--(IV) from Section~3 and \thetag{3-1}.
We also note that $\# F_n=h_n$ and
$$
F_{n-1}+C_{n}-F_{n-1}-C_{n}\subset F_{n}.\tag4-1
$$
Let $(X,\mu,(T_n)_{n\in\Bbb Z})$ denote the dynamical system associated with $(C_n,F_{n-1})_{n\ge 1}$.

\definition{Definition 4.1}
The transformation $T_1$ (or, rather, the dynamical system $(X,\mu, T_1)$)
is called {\it the infinite Chacon transformation}.
\enddefinition

For a (more common) cutting-and-stacking definition of the infinite Chacon transformation we refer to \cite{AdFrSi}, \cite{JaRoRu} and Section~7 of the present paper.

We will need the following simple lemma.

\proclaim{Lemma 4.2}
If $x\in X_{n-1}\cap[f]_n$  and $T_gx\in X_{n-1}$ for  some $g,f\in F_{n}$,
then $g+f\in F_n$.
\endproclaim
\demo{Proof}
Indeed, consider  expansions
$x=(f_{n-1},c_n, c_{n+1},\dots)\in X_{n-1}$
and $T_gx=(f_{n-1}',c_{n}',c_{n+1}',\dots)\in X_{n-1}$ of $x$ and $T_gx$
 respectively.
Then there is $m\ge n$ such that $g=(f'_{n-1}-f_{n-1})+\sum_{l=n}^m(c_l'-c_l)$.
 Moreover, if $m>n$ we can additionally assume that     $c_m\ne c_{m}'$.
It follows from (III) and~\thetag{4-1} that
$g\in F_{m-1}+c_m'-c_m$.
Since $g\in F_n\subset F_{m-1}$, we get a contradiction with (IV).
Hence $m=n$.
Since $x\in[f]_n$, it follows that  $f=f_{n-1}+c_n$.
This implies that $g=f_{n-1}'+c_n-f\in F_n-f$, as desired.
\qed
\enddemo

We first describe the recurrent points of the Cartesian square of $T_1$.

\proclaim{Proposition 4.3}
Let
 $(x,x')\in X_n\times X_n$ for some $n>0$.
Expand $x$ and $x'$ as $x=(f_n,c_{n+1},\dots)$  and  $x'=(f_n',c_{n+1}',\dots)$.
\roster
\item"\rom{(i)}"
If there is $L>n$ such that $\{c_l,c_l'\}=\{-h_{l-1},h_{l-1}+1\}$ for each $l>L$ then the $(T_1\times T_1)$-orbit of $(x,x')$ intersects each compact subset of $X\times X$ at most
 finitely many times.
 \item"\rom{(ii)}"
 Otherwise $(x,x')$ is $(T_1\times T_1)$-recurrent.
 \endroster
 \endproclaim
\demo{Proof}
(i) Let $J:=\{1,\dots,L\}$.
We now let $\alpha_l:=c_l$ and $\beta_l:=c_l'$ for all $n\not\in J$.
It is easy to see that \thetag{3-4} holds.
Then by Theorem~3.2(i), the subset $\psi_{\alpha,\beta}(X^\alpha)$ is closed in $X\times X$. 
However this subset is a single $(T_1\times T_1)$-orbit and the point $(x,x')$ belongs to this orbit.
This follows directly from the definition of  $\psi_{\alpha,\beta}$ (see~\S 3).

(ii)
Let $J:=\{l>0\mid \{c_l,c_l'\}\ne\{-h_{l-1},h_{h-1}+1\}\}$.
We set $\alpha_l:=c_l$ and $
\beta_l:=c_l'$ for each $l\not\in J$.
Since $J$ is infinite, $X^\alpha$ is a perfect space.
It follows now from Theorem~3.2(i) that  $\psi_{\alpha,\beta}$ is a perfect subset of  $X\times X$.
Of course, $(x,x')\in \psi_{\alpha,\beta}(X^\alpha)$.
In view of Theorem~3.2(ii), the $(T_1\times T_1)$-orbit of $(x,x')$ is dense in $\psi_{\alpha,\beta}(X^\alpha)$.
Hence $(x,x')$ is $(T_1\times T_1)$-recurrent.
\qed
\enddemo

\remark{Remark 4.4} In a similar way one can prove the following generalization of Proposition~4.3(i).
Let $d>1$ and let $(x^{(1)},\dots,x^{(d)})\in X^d$.
If there is $L>0$ such that $\{c^{(1)}_l,\dots,c^{(d)}_l\}\supset\{-h_{l-1},h_{l-1}+1\}$ for all $l>L$
then  the $(T_1)^{\times d}$-orbit of $(x^{(1)},\dots,x^{(d)})$ intersects each compact subset of $X^d$ at most
 finitely many times.
 Hence the point $(x^{(1)},\dots,x^{(d)})$ is not $T_1^{\times d}$-recurrent.
 All the other points of $X^d$ are $T_1^{\times d}$-recurrent.

\endremark

\proclaim{Corollary 4.5} The counting measure\footnote{ This measure is infinite, $\sigma$-finite, $(T_1\times T_1)$-invariant, ergodic but non-conservative.} on the $(T_1\times T_1)$-orbit of a point $(x,x')\in X\times X$ is a Radon measure on $X\times X$ if and only if the condition of Proposition~4.3(i) is satisfied.
\endproclaim

We also show that $T_1$ is totally ergodic.

\proclaim{Proposition 4.6} If $m>1$ then the transformation $T_m$ is ergodic.
\endproclaim
\demo{Proof} Let $k>0$.
We set $h:=h_k+h_{k+1}+\cdots +h_{k+m-1}$.
Take $f,f'\in F_k$.
There is $j\in\{0,1,\dots,m-1\}$ such that $f'-f+h+j$ is divisible by $m$.
Take $f,f'\in F_k$
Then $[f-h_{k+j}-
\cdots-h_{k+m-1}]_{k+m}\subset[f]_k$,
$$
\align
T_{h+j+f'-f}[f-h_{k+j}-
\cdots-h_{k+m-1}]_{k+m} &=[f'+(h_k+1)+\cdots+(h_{k+j-1}+1)]_{k+m}\\
&\subset [f']_k
\endalign
$$
and $\mu([f-h_{k+j}-
\cdots-h_{k+m-1}]_{k+m})=3^{-m}\mu([f]_k)$.
It remains to apply Lemma~3.1.
\qed
\enddemo

\head 5. Radon invariant measures  for~Cartesian~powers~of the infinite Chacon transformation
\endhead

We will  need the following general lemma about splitting of Radon  measures into direct products (cf. with \cite{RuSi, Lemma~3.1.1}).
It is, in fact, equivalent to \cite{JaRoRu, Lemma~A.1}.
We provide  a short proof of it.

\proclaim{Lemma 5.1}  
Let  $T$ and $S$ be homeomorphisms of locally compact Cantor spaces $X$ and $Y$ respectively.
If $\lambda\in\Cal M^e_{\text{\rom{Ra}}}(X\times Y,T\times S)$ and $\lambda\circ
(\text{\rom{Id}}\times S)=\lambda$
then there exist $\mu\in\Cal M^e_{\text{\rom{Ra}}}(X,T)$ and
$\nu\in \Cal M^e_{\text{\rom{Ra}}}(Y,S)$ such that $\lambda=\mu\times\nu$.
\endproclaim
\demo{Proof}
Let $\widetilde\mu$ be a probability measure on $(X,\goth B)$ which is equivalent (i.e. it has the same ideal of Borel subsets of $0$-measure) to the projection of $\lambda$ to  $X$\footnote{Such a measure can be obtained as the projection (to $X$)  of a probability measure on $X\times Y$ equivalent to $\lambda$.}.
Then $\widetilde\mu$ is quasi-invariant and ergodic under $T$.
Consider the disintegration of $\lambda$ with respect to $\widetilde\mu$:
$$
\lambda=\int_X\delta_x\times\lambda_x\,d\widetilde\mu(x),
$$
where the map $X\ni x\mapsto\lambda_x$ is the corresponding system of $\sigma$-finite measures on $Y$.
Since $\lambda\in \Cal M_{\text{\rom{Ra}}}(X\times Y)$, it follows that $\lambda_x\in \Cal M_{\text{\rom{Ra}}}(X) $ for $\widetilde\mu$-a.e. $x$.
We have that
$$
\lambda\circ(T\times S)=\int_X\delta_{T^{-1}x}\times\lambda_x\circ S\,d\widetilde\mu(x)=\int_X\delta_x\times
\frac{d\widetilde\mu\circ T}{d\widetilde\mu}(x)\cdot 
\lambda_{Tx}\circ S\,d\widetilde\mu(x).
$$
In view of the uniqueness of the disintegration we obtain that $\lambda$ is invariant under $T\times S$ if and only if $\lambda_x=\frac{d\widetilde\mu\circ T}{d\widetilde\mu}(x)\cdot \lambda_{Tx}\circ S$ for a.e. $x$.
In a similar way, $\lambda$ is invariant under Id$\times S$ if and only if $\lambda_x=\lambda_x\circ S$ for a.e. $x$.
Hence $\lambda_x=\frac{d\widetilde\mu\circ T}{d\widetilde\mu}(x)\cdot \lambda_{Tx}$ for a.e. $x$.
Since $T$ is ergodic and $\lambda_x$ is a Radon measure for a.e. $x\in X$,
it follows that there is a single Radon measure $\nu$ on $Y$ and a measurable map $X\ni x \mapsto a(x)\in \Bbb R^*_+$ such that $\lambda_x=a(x)\nu$ for $\mu$-a.e. $x$.
We now obtain that $\nu \circ S =\nu$ and $a(x)=\frac{d\widetilde\mu\circ T}{d\widetilde\mu}(x) a(Tx)$ for a.e. $x$.
Define a $\sigma$-finite measure measure $\mu$ on $X$ by setting $\frac{d\mu}{d\widetilde\mu}(x)=a(x)$ for all $x\in X$.
Then $\mu$ is  invariant under $T$.
 Moreover, $\lambda=\mu\times\nu$ and hence $\mu\in \Cal M_{\text{Ra}}^e(X,T)$ and 
 $\nu\in \Cal M_{\text{Ra}}^e(Y,S)$.
\qed
\enddemo

We now state one of the main results of the paper.

\proclaim{Theorem 5.2}
The infinite Chacon transformation
 has  Radon MSJ.
 \endproclaim

\demo{Proof}
Let $(X,\mu, T_1)$ denote the infinite Chacon transformation as above.
 We fix $d>1$ and take $\lambda\in \Cal M^{e}_{\text{Ra}}(X^d,(T_1)^{\times d})$.
Let  $z=(x^{(1)},\dots,x^{(d)})\in \Cal G(\lambda)$.
We recall that $X=\bigcup_{n\ge 0}X_n$ and $X_n:=F_n\times C_{n+1}\times\cdots$.
If $m$ is large so that $z\in (X_m)^d$, we consider expansion
$x^{(j)}=(f_m^{(j)},c_{m+1}^{(j)},c_{m+2}^{(j)},\dots)\in X_m$ of $x^{(j)}$ for each $j=1,\dots,d$.
We call $z$ {\it extreme} if for each sufficiently large $n>m$,
either $\{-h_{n-1},h_{n-1}+1\}\subset \{c_n^{(1)},\dots,c_n^{(d)}\}$
or $c_n^{(1)}=\cdots=c_n^{(d)}$.

We first show that if $z$ is not extreme then $\lambda$ splits into direct product of two its marginals.
Indeed, there is
 an infinite subset $\Cal N\subset\Bbb N$ and a non-empty  subset $J\subset \{1,\dots,d\}$ such that $c^{(j)}_n=0$ if $j\in J$ and either $c^{(j)}_n=-h_{n-1}$ for all $j\not\in J$ or $c^{(j)}_n=h_{n-1}+1$ for all $j\in\{1,\dots,d\}\setminus J$ for each $n\in\Cal N$.
We consider only the former case since the latter one is similar.
Let $A^{(1)},\dots,A^{(d)}$ be compact open subsets in $X$.
Then there is $r>0$ such that these subsets are $r$-cylinders.
Hence for each  $n\ge r$ and each $j=1,\dots,d$, there is a subset $A^{(j)}_n\subset F_n$ such that $A^{(j)}=[A^{(j)}_n]_n$.
It follows from Lemma~4.2 that
$$
\sum_{i\in F_n}1_{A^{(1)}\times\cdots\times A^{(d)}}((T_1^{\times d})^iz)=\#\bigg(F_n\cap\bigcap_{j=1}^d(A_n^{(j)}-f_n^{(j)})\bigg).\tag5-1
$$
Since we may assume without loss of generality that $A^{(j)}\cup\{x^{(j)}\}\subset X_{n-1}$, it follows from \thetag{4-1} that
$$
A_n^{(j)}-f_n^{(j)}\subset F_n\quad\text{for each $j=1,\dots,d.$}\tag5-2
$$
 By the definition of generic point,  for each $\epsilon>0$ and each sufficiently large $n$, we have (in view of \thetag{5-1} and \thetag{5-2}) that
$$
\frac{\#\bigg(\bigcap_{j=1}^d(A_n^{(j)}-f_n^{(j)})\bigg)}{\#\bigg(\bigcap_{j=1}^d
(F_{r,n}-f_n^{(j)})\bigg)}=\frac{\lambda(A^{(1)}\times\cdots\times A^{(d)})}{\lambda(X_r^d)}\pm \epsilon\tag5-3
$$
where $F_{r,n}:=F_r+C_{r+1}+\cdots+C_n$.
Now choose $n$ so that  $n+1\in\Cal N$.
Then
$$
f_{n+1}^{(j)}=f_{n}^{(j)}+c_{n+1}^{(j)}=
\cases
f_{n}^{(j)}, &\text{if }j\in J \\
f_{n}^{(j)}-h_n ,&\text{if }j\notin J .
\endcases
$$
Since  $A_{n+1}^{(j)}=A_n^{(j)}+C_{n+1}$,
we obtain that
$$
\align
\bigcap_{j=1}^d(A_{n+1}^{(j)}&-f_{n+1}^{(j)})=\bigcap_{j\in J}(A_{n}^{(j)}-f_{n}^{(j)}+C_{n+1})\cap\bigcap_{j\not\in J}(A_{n}^{(j)}-f_{n}^{(j)}+C_{n+1}+h_n)\\
&=\bigcap_{j\in J}(A_{n}^{(j)}-f_{n}^{(j)}+\{0,h_n+1\})\cap\bigcap_{j\not\in J}(A_{n}^{(j)}-f_{n}^{(j)}+\{0,h_n\})\\
&=\bigcap_{j=1}^d(A_{n}^{(j)}-f_{n}^{(j)})\sqcup\left(h_n+
\bigcap_{j\in J}(A_{n}^{(j)}-f_{n}^{(j)}+1)\cap\bigcap_{j\not\in J}(A_{n}^{(j)}-f_{n}^{(j)})\right).
\endalign
$$
We define a map  $S:X^d\to X^d$ by setting
$S(x^{(1)},\dots,x^{(d)})=(y^{(1)},\dots,y^{(d)})$, where $y^{(j)}:=x^{(j)}$ if $j\in J$ and
$y^{(j)}:=T_1x^{(j)}$ if $j\not\in J$.
Of course, $S$ is a homeomorphism of $X^d$.
It commutes   with $T_1^{\times d}$.
Let
$$
g^{(j)}:=\cases
0 , &\text{if }j\in J \\
1, &\text{if }j\notin J .
\endcases
$$
Then
$$
\#\bigg(\bigcap_{j=1}^d(A_{n+1}^{(j)}-f_{n+1}^{(j)})\bigg)=
\#\bigg(\bigcap_{j=1}^d(A_{n}^{(j)}-f_{n}^{(j)})\bigg)+
\#\bigg(\bigcap_{j=1}^d(A_{n}^{(j)}+g^{(j)}-f_{n}^{(j)})\bigg).
$$
In a similar way we obtain that
$$
\#\bigg(\bigcap_{j=1}^d
(F_{r,n+1}-f_{n+1}^{(j)})\bigg)=
\#\bigg(\bigcap_{j=1}^d
(F_{r,n}-f_n^{(j)})\bigg)+
\#\bigg(\bigcap_{j=1}^d
(F_{r,n}+g^{(j)}-f_n^{(j)})\bigg).
$$
Thus we obtain that the lefthand side of \thetag{5-3} (with $n+1$ in place of $n$) equals
$$
\frac{\#\bigg(\bigcap_{j=1}^d(A_{n}^{(j)}+g^{(j)}-f_n^{(j)})\bigg)
+\#\bigg(\bigcap_{j=1}^d(A_{n}^{(j)}-f_n^{(j)})\bigg)}
{\#\bigg(\bigcap_{j=1}^d
(F_{r,n}+g^{(j)}-f_n^{(j)})\bigg)+\#\bigg(\bigcap_{j=1}^d
(F_{r,n}-f_{n}^{(j)})\bigg)}.\tag5-4
$$
It is straightforward to verify that   $\lambda(S(X_r)^d)/\lambda((X_r)^d)=1\pm\epsilon$ if $r$ is large enough.
  Therefore it follows from \thetag{5-3}  and \thetag{5-4} that
$$
\frac{\lambda(A^{(1)}\times\cdots\times A^{(d)})}{\lambda(X_r^d)}\pm \epsilon=
\frac{\frac{\lambda\circ S(A^{(1)}\times\cdots\times A^{(d)})}{\lambda(X_r^d)}+\frac{\lambda(A^{(1)}\times\cdots\times A^{(d)})}{\lambda(X_r^d)}\pm2\epsilon}{2\pm\epsilon}.
$$
Hence $\lambda=(\lambda+\lambda\circ S)/2$, i.e. $\lambda\circ S=\lambda$.
Lemma~5.1 yields that $\lambda=\lambda_1\times\lambda_2$, where 
$\lambda_1\in\Cal M_{\text{Ra}}^e(X^{J},T_1^J)$ and
$\lambda_2\in\Cal M_{\text{Ra}}^e(X^{{\{1,\dots,d\}\setminus J}},T_1^{{\{1,\dots,d\}\setminus J}})$.
Continuing this way several times we  obtain finally a splitting  of $\lambda$ into
direct product of  its marginals  whose generic points are all extreme.

Thus to complete the proof of the thorem it suffices to prove the following fact:
if every $\lambda$-generic point is extreme then either
there is $s\in\{1,\dots,d\}$, $d>1$, such that the projection of $\Cal G(\lambda)$ to the $s$-coordinate is $\mu$-negligible or
there exist $n_1,\dots,n_{d-1}\in\Bbb Z$ such that $\lambda=\mu_{T^{n_1},\dots, T^{n_{d-1}}}$\footnote{We omit the case where $d=1$ because it is trivial. We recall that every $(C,F)$-action is Radon strictly ergodic.}.
Thus we fix $z\in\Cal G(\lambda)$ as above and 
find $m>0$ such that for each $n>m$, either $\{-h_{n-1}, h_{n-1}+1\}\subset \{c_n^{(1)},\dots, c_n^{(d)}\}$ or
$c_n^{(1)}=\dots= c_n^{(d)}$.
Put
$$
I:=\{i>m\mid c^{(1)}_i=\cdots=c^{(d)}_i\}.
$$
We can now describe explicitly the set $\Cal G(\lambda)$.
Namely, we claim that
$$
\aligned
\Cal G(\lambda)&\cap([f_m^{(1)}]_m  \times \cdots\times[f_m^{(d)}]_m)\\
&=\bigcap_{l>m}\bigcap_{j_1,j_2=1}^d\{\widetilde z\in (X_m)^d\mid \widetilde c^{(j_1)}_l-\widetilde c^{(j_2)}_l 
=
c^{(j_1)}_l- c^{(j_2)}_l\},
\endaligned
\tag5-5
$$
where the coordinates $\widetilde c^{(j)}_i$ are taken from the expansion
$$
\widetilde z=((\widetilde f_m^{(1)},\widetilde c_{m+1}^{(1)},\widetilde c_{m+2}^{(1)},\dots),\dots, (\widetilde f_m^{(d)},\widetilde c_{m+1}^{(d)},\widetilde c_{m+2}^{(d)},\dots))
$$
 of $\widetilde z$.
Indeed, let
 $A^{(1)},\dots, A^{(d)}$ be  arbitrary $r$-cylinders with $r\ge m$.
Then, as above, we have
$$
\frac{\sum_{i\in F_{r+1}}1_{A^{(1)}\times\cdots\times A^{(d)}}((T^{\times d})^iz)}
{\sum_{i\in F_{r+1}}1_{X_r^d}((T^{\times d})^iz)}
=
\frac{
\#\bigg(\bigcap_{j=1}^d(A_{r+1}^{(j)}-f_{r+1}^{(j)})\bigg)}
{\#\bigg(\bigcap_{j=1}^d(F_{r,{r+1}}-f_{r+1}^{(j)})\bigg)}.\tag5-6
$$
Consider now 2 cases.
If $r+2\in I$ then
$$
\bigcap_{j=1}^d(A_{r+2}^{(j)}-f_{r+2}^{(j)})=
\bigcap_{j=1}^d(A_{r+1}^{(j)}-f_{r+1}^{(j)}+C_{r+2}-c_{r+2}^{(j)})=
\bigsqcup_{c\in C_{r+2}-c_{r+2}^{(1)}}\bigcap_{j=1}^d(A_{r+2}^{(j)}-f_{r+2}^{(j)}+c)
$$
and hence $\#(\bigcap_{j=1}^d(A_{r+2}^{(j)}-f_{r+2}^{(j)}))
=3\#(\bigcap_{j=1}^d(A_{r+1}^{(j)}-f_{r+1}^{(j)}))$.
In a similar way, $\#(\bigcap_{j=1}^d(F_{r,{r+2}}-f_{r+2}^{(j)}))=
3\#(\bigcap_{j=1}^d(F_{r,{r+1}}-f_{r+1}^{(j)}))$.

If $r+2\not\in I$ then
$$
\bigcap_{j=1}^d(A_{r+2}^{(j)}-f_{r+2}^{(j)})=
\bigcap_{j=1}^d(A_{r+1}^{(j)}-f_{r+1}^{(j)}+C_{r+2}-c_{r+2}^{(j)})=
\bigcap_{j=1}^d(A_{r+1}^{(j)}-f_{r+1}^{(j)})
$$
and hence $\#(\bigcap_{j=1}^d(A_{r+2}^{(j)}-f_{r+2}^{(j)}))
=\#(\bigcap_{j=1}^d(A_{r+1}^{(j)}-f_{r+1}^{(j)}))$.
In a similar way, $\#(\bigcap_{j=1}^d(F_{r,{r+2}}-f_{r+2}^{(j)}))=
\#(\bigcap_{j=1}^d(F_{r,{r+1}}-f_{r+1}^{(j)}))$.
Thus in every of the two cases we obtain that
$$
\frac{
\#\bigg(\bigcap_{j=1}^d(A_{r+1}^{(j)}-f_{r+1}^{(j)})\bigg)}
{\#\bigg(\bigcap_{j=1}^d(F_{r,{r+1}}-f_{r+1}^{(j)})\bigg)}=
\frac{
\#\bigg(\bigcap_{j=1}^d(A_{r+2}^{(j)}-f_{r+2}^{(j)})\bigg)}
{\#\bigg(\bigcap_{j=1}^d(F_{r,{r+2}}-f_{r+2}^{(j)})\bigg)}
\tag5-7
$$
for each $r\ge m$.
Moreover, 
if we  assume additionally that $A^{(j)}\subset X_{r-1}$ or, equivalently, 
$A_r^{(j)}\subset F_{r-1}+C_r$  for each $j=1,\dots,d$ then the same argument 
as above yields that \thetag{5-7} holds also 
if we replace $r$ with $r-1$.
It follows from this, \thetag{5-6} and the fact that $z$ is $\lambda$-generic that
$$
\frac{\lambda(A^{(1)}\times\cdots\times A^{(d)})}{\lambda(X_r^d)}=
\frac{
\#\bigg(\bigcap_{j=1}^d(A_{r}^{(j)}-f_{r}^{(j)})\bigg)}
{\#\bigg(\bigcap_{j=1}^d(F_{r}-f_{r}^{(j)})\bigg)}.
\tag5-8
$$
If  $\widetilde z$ is another $\lambda$-generic point 
with $\widetilde z\in[f_m^{(1)}]_m  \times \cdots\times[f_m^{(d)}]_m$
then it is extreme and therefore  \thetag{5-8} holds with $\widetilde f_{r}^{(j)}$ in place  of  $f_{r}^{(j)}$, $j=1,\dots,d$.
This yields
$$
\frac{
\#\bigg(\bigcap_{j=1}^d(A_{r}^{(j)}-f_{r}^{(j)})\bigg)}
{\#\bigg(\bigcap_{j=1}^d(F_{r}-f_{r}^{(j)})\bigg)}=
\frac{
\#\bigg(\bigcap_{j=1}^d(A_{r}^{(j)}-\widetilde f_{r}^{(j)})\bigg)}
{\#\bigg(\bigcap_{j=1}^d(F_{r}-\widetilde f_{r}^{(j)})\bigg)}.
\tag5-9
$$
We obtain that
$\bigcap_{j=1}^d(A^{(j)}-f_{r}^{(j)})\ne\emptyset$ if and only if
$\bigcap_{j=1}^d(A^{(j)}-\widetilde f_{r}^{(j)})\ne\emptyset$.
Substituting $A^{(j)}:=\{f_r^{(j)}\}$, we obtain that
$f_r^{(j_1)}-f_r^{(j_2)}=\widetilde f_r^{(j_1)}-\widetilde f_r^{(j_2)}$ for each $1\le j_1<j_2\le d$.
Since $r$ is arbitrary, we conclude that
$c_l^{(j_1)}-c_l^{(j_2)}=\widetilde c_l^{(j_1)}-\widetilde c_l^{(j_2)}$ for each $1\le j_1<j_2\le d$ and all sufficiently large $l$.
This proves the  ``$\subset$''-part of the equality in \thetag{5-5}.
The ``$\supset$''-part of this equality follows easily from \thetag{5-9} and~\thetag{5-6}.
Thus \thetag{5-5} is proved.

Consider now two cases.
If $\Bbb N\setminus I$ is finite then
 we have $x^{(j)}=T_{g_j}x^{(1)}$ for $g_j:=f_m^{(j)}-f_m^{(1)}+\sum_{n>m}(c^{(j)}_n-c^{(1)}_n)$, $j=2,\dots,d$.
Hence for each sufficiently large $r$ and each $j\in\{1,\dots,d\}$, we have that $f^{(j)}_r=g_j+f_r^{(1)}$.
Then \thetag{5-7} implies that
$$
\frac{\lambda(A^{(1)}\times\cdots\times A^{(d)})}{\lambda(X_r^d)}=
\frac{
\#\bigg(\bigcap_{j=1}^d(A_{r}^{(j)}-g_j)\bigg)}
{\#\bigg(\bigcap_{j=1}^d(F_{r}-g_j)\bigg)}=
\frac{\mu_{T_{g_2},\dots, T_{g_d}}(A^{(1)}\times\cdots\times A^{(d)})}
{\mu_{T_{g_2},\dots, T_{g_d}}(X_r^d)}.
$$
Hence $\lambda$ is a multiple of $\mu_{T_{g_2},\dots, T_{g_d}}$.
Thus $\lambda$ is a graph-joining.
Conversely, if $\lambda$ is a graph-joining then $\Bbb N\setminus I$ is finite.

Now consider the second  case, where $\Bbb N
\setminus I$ is infinite.
Then it follows from \thetag{5-5} that
 there exists $s\in\{1,\dots,d\}$ and an infinite subset $I_0\subset\{m+1,m+2,\dots\}\setminus I$ such that for  every $\widetilde z\in\Cal G(\lambda)\cap ([f_m^{(1)}]_m  \times \cdots\times[f_m^{(d)}]_m)$, we have that
$\widetilde c_j^{(s)}\ne 0$ whenever $j\in I_0$.
However,  $\mu(\{x=(f_m,c_{m+1},\dots)\in X_m\mid  c_j\ne 0\text{ for all $j\in I_0$}\})=0$.
Thus the projection of the subset $\Cal G(\lambda)\cap ([f_m^{(1)}]_m  \times \cdots\times[f_m^{(d)}]_m)\subset X^d$ to the $s$-th coordinate is
$\mu$-negligible.
Since the projection is equivariant and $\lambda([f_m^{(1)}]_m  \times \cdots\times[f_m^{(d)}]_m)>0$, it follows that the projection of the entire $\Cal G(\lambda)$
to the $s$-th coordinate is also
$\mu$-negligible.
Hence  $\lambda\not\in J_{d,\text{Ra}}((T_1)^{\times d})$, as desired.
Thus, $T_1$ has Radon MSJ.
\qed

\enddemo

 We note that  not only the Radon MSJ property for $T_1$ was established  in the proof of Theorem~5.2  but also
 a complete description of the set of ergodic $(T_1)^{\times d}$-invariant Radon measures on $X^d$ was, in fact, obtained.
 Namely, the following theorem was proved indeed.
 
 \proclaim{Theorem 5.3} Let $\lambda\in \Cal M_{\text{\rom{Ra}}}^e(X^d,(T_1)^{\times d})$.
 Then there is a partition of $\{1,\dots,d\}$ into subsets $J_1,\dots, J_k$ such that $\lambda$ splits into direct product of its marginals $\lambda_i$ on $X^{J_i}$, each $\lambda_i$ is (up to multiplicative constant) either $\mu$ if $\# J_i=1$  or an invariant quasi-graph Radon measure
 if $\# J_i>1$.
 \endproclaim

We now show that there are non-ergodic 2-fold Radon self-joinings $\rho$ of $T_1$ such that almost every ergodic component of $\rho$ is not a joining of $T$:
the ergodic components  of $\rho$ are conservative quasi-graph invariant Radon measures which are not graph-joinings.

\example{Example 5.4}
Let $\boldsymbol{A}:=\{-1,0,1\}^\Bbb N$.
Given a sequence $a:=(a_n)_{n=1}^\infty\in\boldsymbol{A}$, we define a measure $\lambda_a$
on $X_0^2$ considered as the infinite product $\{(0,0)\}\times C_1^2\times C_2^2\times\cdots$
by setting
$$
\lambda_a=\delta_{(0,0)}\otimes\bigotimes_{n>0}\gamma_n^a,\tag5-10
$$
where $\gamma_n$ is a probability on $C_n^2$ given by
$$
 \gamma_n^a:=
 \cases
 \delta_{(-h_{n-1},h_{n-1}+1)} &\text{if }a_n=-1, \\
 \delta_{(h_{n-1}+1,-h_{n-1})} &\text{if }a_n=1, \\
 \frac 13\sum_{c\in C_n}\delta_{(c,c)} &\text{if }\alpha_n=0.
 \endcases
 \tag5-11
 $$
We recall that $C_n=\{-h_{n-1},0,h_{n-1}+1\}$.
Let $J_a:=\{n\mid a_n\ne 0\}$.
We now set for each $n\not\in J_a$,
$$
\alpha_n^a:=
\cases
-h_{n-1},&\text{if }a_n=-1\\
h_{n-1}+1,&\text{if }a_n=1
\endcases,
\quad
\beta_n^a:=
\cases
h_{n-1}+1,&\text{if }a_n=-1\\
-h_{n-1},&\text{if }a_n=1
\endcases.
$$
It is easy to see that the sequences $\alpha^a:=(\alpha^a_n)_{n\not\in J}$ and $\beta^a:=(\beta_n^a)_{n\not\in J}$ satisfy~\thetag{3-4}.
Then it follows from Theorem~3.2(iv) that $\lambda_a=\mu^{\alpha^a}\circ\psi^{-1}_{\alpha^a,\beta^a}\restriction X_0^2$ and that $\mu^{\alpha^a}\circ\psi^{-1}_{\alpha^a,\beta^a}$ is a $(T_1\times T_1)$-invariant  quasi-graph Radon measure\footnote{For the definition of $\mu^\alpha$ and $\psi_{\alpha,\beta}$ we refer to \S 3.2.}.
For each $n>0$, let  $\kappa_n$ be a probability on $\{-1,0,1\}$ such that
$\kappa_n(-1)=\kappa_n(1)=\frac 1{2n}$ and $\kappa_n(0)=1-\frac 1n$.
We now set $\kappa:=\bigotimes_{n>0}\kappa_n$ and $\lambda:=\int_{\boldsymbol A}\lambda_a\,d\kappa(a)$.
Then 
$$
\lambda=\left(\int_{\boldsymbol A}\mu^{\alpha^a}\circ\psi^{-1}_{\alpha^a,\beta^a}\,d\kappa(a)\right)\restriction X_0^2.
$$
We claim that the Radon measure 
$$
\rho:=\int_{\boldsymbol A}\mu^{\alpha^a}\circ\psi^{-1}_{\alpha^a,\beta^a}\,d\kappa(a)\tag5-12
$$
is a 2-fold self-joining of $T_1$.
For that we have to verify that the two coordinate projections of this measure are equivalent to $\mu$.
It suffices to show that the two coordinate projections of $\lambda$ are equivalent to $\mu\restriction X_0$.
A straightforward computation shows that the projection of $\lambda$ to the first coordinate is
the following measure\footnote{We use \thetag{5-10} and \thetag{5-11} to obtain this.}:
$$
\delta_0\otimes\bigotimes_{n>0}\left(\kappa_n(-1)\delta_{-h_{n-1}}+\kappa_n(1)\delta_{h_{n-1}+1}+\frac{\kappa_n(0)}3\sum_{c\in C_n}\delta_c \right).
$$
 It is equivalent to the measure $\delta_0\otimes\bigotimes_{n>0}\left(\frac{1}3\sum_{c\in C_n}\delta_c \right)$ (i.e. to $\mu\restriction X_0$)
 by the Kakutani theorem on equivalence of infinite product measures \cite{Ka}.
 In a similar way one can prove that the projection of $\rho$ to the second coordinate is  also equivalent to $\mu$.
 Therefore $\rho\in J_{2,\text{Ra}}(T_1)$.
On the other hand, 
$$
\align
0&=\kappa(\{a\mid J_a\text{ is cofinite}\})=\kappa(\{a\mid \mu^{\alpha^a}\circ\psi^{-1}_{\alpha^a,\beta^a}\text{ is a 2-fold joining of $T_1$}\}),\\
0&=\kappa(\{a\mid J_a\text{ is finite}\})=\kappa(\{a\mid \mu^{\alpha^a}\circ\psi^{-1}_{\alpha^a,\beta^a}\text{ is purely atomic}\}).
\endalign
$$
We note that \thetag{5-12} is the ergodic decomposition of $\rho$ and $\kappa$ is the corresponding measure on the space of ergodic components of $\rho$.
Hence almost all ergodic components of $\rho$ are conservative Radon dynamical systems whose coordinate projections are singular to $\mu$.
\endexample

\head 6. Uncountable family of pairwise~Radon~disjoint infinite  Chacon like transfrmations
\endhead

Take an infinite sequence $\omega\in\{0,1\}^\Bbb N$.
Define a sequence of positive integers $(h_n)_{n=0}^\infty$ by setting $h_0:=1$ and $h_{n+1}:=6h_n+1$ for $n>0$.
For each $n\ge 0$, we now let
$$
\align
F_n&:=\left\{-\frac{h_n-1}2,-\frac{h_n-3}2,\dots,\frac{h_n-1}2\right\}\qquad\text{ and }\\
C_{n+1}^\omega&:=\{-h_n-1+\omega(n),0,h_n+\omega(n)\}.
\endalign
$$
The sequence $(C_n^\omega,F_{n-1})_{n\ge 1}$ satisfies \thetag{1-1} and (I)--(IV) from Section~1.
We also note that $\# F_n=h_n$
and
$$
F_n+C_{n+1}^\omega-F_n-C_{n+1}^\omega\subset F_{n+1}.
$$
Let $(X^\omega,\mu^\omega,T^\omega)$ denote the associated $(C,F)$-dynamical system.
It is Radon strictly ergodic.
We call it an {\it infinite Chacon like} transformation.
The infinite Chacon transformation corresponds to the case $\omega(n)=1$ for all $n$.

If $\omega,\omega'\in\{0,1\}^\Bbb N$, we write $\omega\sim\omega'$ if the pair $(\omega,\omega')$ belongs to the tail equivalence relation on $\{0,1\}^\Bbb N$, i.e. there is $N>0$ such that $\omega(n)=\omega'(n)$ for all $n>N$.

\proclaim{Theorem 6.1}
\roster
\item"\rom(i)"
For each $\omega\in \{0,1\}^\Bbb N$, the  dynamical system $(X^\omega,\mu^\omega,T^\omega_1)$ is totally ergodic.
 It  has  Radon MSJ.
$C(T_1^\omega)=\{T_n^\omega\mid n\in\Bbb Z\}$.
\item"\rom(ii)"
If $\omega\sim\omega'$ then  $(X^\omega,\mu^\omega,T^\omega_1)$ is isomorphic to
 $(X^{\omega'},\mu^{\omega'},T^{\omega'}_1)$.
 \item"\rom(iii)"
If $\omega\not\sim\omega'$ then  $(X^\omega,\mu^\omega,T^\omega_1)$ is Radon disjoint with
 $(X^{\omega'},\mu^{\omega'},T^{\omega'}_1)$.
\endroster
\endproclaim

\demo{Sketch of the proof}
(i) is proved in the same way as Theorem~5.2 and Proposition~4.6.

(ii) Let $\omega(i)=\omega'(i)$ for all $i>N$.
Given $x\in X^\omega$, let $x=(f_n,c_{n+1},\cdots)\in (X^\omega)_n$ for some  $n>N$.
We now define $\phi:X^\omega\to X^{\omega'}$ by setting
$\phi(x):=(f_n,c_{n+1},\dots)\in  (X^{\omega'})_n$.
Then $\phi$ is an isomorphism  of $(X^\omega,\mu^\omega, T^\omega_1)$ onto $(X^{\omega'},\mu^{\omega'},T^{\omega'}_1)$.

(iii) Let $\lambda$ be an ergodic $(T^\omega_1\times T^{\omega'}_1)$-invariant measure on $X^\omega\times X^{\omega'}$ whose marginals are equivalent to $\mu^\omega$ and $\mu^{\omega'}$
respectively.
Fix  a generic point $(x,x')$ of $\lambda$.
Find $m$ such that $x\in (X^\omega)_m$ and $x'\in (X^{\omega'})_m$.
Consider expansions $x=(f_m,c_{m+1},\dots)$ and $x'=(f_m',c_{m+1}',\dots)$ of $x$ and $x'$ in $(X^\omega)_m$ and $(X^\omega)_m$ respectively.
Let
$$
J:=\{i\in\Bbb N\mid \omega(i)\ne\omega'(i)\}.
$$
This set is infinite.
One of the following cases takes place.
\roster
\item"(A)"  The subset $J_A:=\{j\in J\mid  \{c_j,c_j'\}=\{h_{j-1}+1,h_{j-1}\}\}$ is infinite.
\item"(B)"  The subset $J_B:=\{j\in J\mid  \{c_j,c_j'\}=\{-h_{j-1},-h_{j-1}-1\}\}$ is infinite.
\item"(C)" The subset $J_C:=\{j\not\in J\mid 0\in\{c_j,c_j'\}\ne\{0\}\}$
is infinite.
\item"(D)"  The subset $J_A\cup J_B\cup J_C$ is finite.
\endroster

We consider separately every case.

Case (A). Let $A$ and $A'$ be two $r$-cylinders with $r>m$.
Then $A=[A_n]_n$ and $A'=[A'_n]_n$ for some subsets $A_n\subset F_n$ and $A_n'\subset F_n'$ for each $n\ge r$.
As in the proof of Theorem~5.2, we obtain that
$$
\sum_{i\in F_n}1_{A\times A'}((T^{\omega}_1\times T^{\omega'}_1)^i(x,x'))=\#((A_n-f_n)\cap(A_n'-f_n')).
$$
Therefore  for each $\epsilon>0$ and each sufficiently large $n$,
we have that
$$
\frac{\#((A_n-f_n)\cap(A_n'-f_n'))}{\#(
(F_{r,n}-f_n)\cap (F_{r,n}-f_n'))}=\frac{\lambda(A\times A')}{\lambda(X_r)}\pm \epsilon.\tag6-1
$$
This follows from the fact that $(x,x')\in\Cal G(T^\omega_1\times T^{\omega'}_1)$.
Now choose $n$ so that  $n+1\in J_A$.
Then  $c_{n+1}=h_{n}+1$ and $c_{n+1}'=h_{n}$ (or  $c_{n+1}=h_{n}$ and $c_{n+1}'=h_{n}+1$, which is considered in a similar way) and hence
 $f_{n+1}=f_n+h_n+1$ and $f_{n+1}'=f_n'+h_n$.
 We now have
$$
(A_{n+1}-f_{n+1})\cap(A_{n+1}'-f_{n+1}')=(A_n-1-f_n+C_{n+1}^\omega-h_n)\cap
(A_n-f_n+C_{n+1}^{\omega'}-h_n).
$$
Since $C^\omega_{n+1}=\{-h_n,0,h_n+1\}$
and $C^{\omega'}_{n+1}=\{-h_n-1,0,h_n\}$\footnote{Because $h_{n}+1\in C^\omega_{n+1}$ and  $h_n\in C_{n+1}^{\omega'}$.}, it follows that
$$
\multline
\#((A_{n+1}-f_{n+1})\cap(A_{n+1}'-f_{n+1}') )
=
\#((A_{n}-f_{n})\cap(A_{n}'-f_{n}') )\\+
\#((A_{n}-1-f_{n})\cap(A_{n}'-f_{n}') )
+
\#((A_{n}-f_{n})\cap(A_{n}'-f_{n}') )
\endmultline
$$
In a similar way we obtain that
$$
\align
\#(
(F_{r,n+1}-f_{n+1})\cap (F_{r,n+1}-f_{n+1}'))
&=2\#(
(F_{r,n}-f_n)\cap (F_{r,n}-f_n'))\\
&+
\#(
(F_{r,n}-1-f_n)\cap (F_{r,n}-f_n')).
\endalign
$$
Applying \thetag{6-1} we now obtain that
$$
\frac{\lambda(A\times A')}{\lambda(X_r)}\pm \epsilon=
\frac{\frac{\lambda(A\times A')}{\lambda(X_r)}+\frac{2\lambda\circ (T^\omega_{-1}\times \text{Id})(A\times A')}{\lambda(X_r)}\pm3\epsilon}{3\pm\epsilon}.
$$
This yields that $\lambda\circ (T^\omega_{-1}\times \text{Id})=\lambda$.
By Lemma 5.1, $\lambda$ equals $\mu\times\mu'$ (up to a multiplicative constant).

 Case (B) is analogous to Case (A).

Case (C) was considered, in fact, in the first part of the proof of Theorem~5.2.
In this case we also obtain that $\lambda =q\cdot\mu\times\mu'$ for some $q>0$.

 Case (D). We will show that in this case one of the two marginals of the set $\Cal G (\lambda)$ is of null measure.
For that we will argue as in the final part of the proof of Theorem~5.2.
Without loss of generality we may assume that
 $\{c_n,c_n'\}\cap (J_A\cup J_B\cup J_C)=\emptyset$ for each $n>m$.
We claim that
$$
\Cal G(\lambda)\cap([f_m]_m\times[f']_m)=
\bigcap_{l>m}\{(\widetilde x,\widetilde x')\in (X^\omega)_m\times (X^{\omega'})_m\mid \widetilde c_l-\widetilde c'_l =
c_l- c'_l\},
\tag6-2
$$
where the coordinates $\widetilde c_l$ and $\widetilde c_l'$ are taken from the expansion
$$
\widetilde x=(\widetilde f_m,\widetilde c_{m+1},\widetilde c_{m+2},\dots), \quad
\widetilde x'=(\widetilde f_m',\widetilde c_{m+1}',\widetilde c_{m+2}',\dots)
$$
 of $\widetilde x$ and $\widetilde x'$.
We let $\alpha_n:=\#((A_{n}-f_{n})\cap(A_{n}'-f_{n}') )$.
One can verify (as in the proof of Theorem~5.2) that
$$
\alpha_{n+1}=
\cases
3\alpha_n&\text{if $c_{n+1}=c_{n+1}'$},\\
\alpha_n&\text{if $0\not\in \{c_{n+1},c_{n+1}'\}$ but $c_{n+1}\ne c_{n+1}'$},\\
2\alpha_n &\text{if  $0\in \{c_{n+1},c_{n+1}'\}$ but $c_{n+1}\ne c_{n+1}'$}.
\endcases
\tag6-3
$$
This implies
$$
\frac{
\#((A_{r+1}-f_{r+1})\cap(A_{r+1}'-f_{r+1}')}
{\#((F_{r,r+1}-f_{r+1})\cap(F_{r,r+1}'-f_{r+1}')}=
\frac{
\#((A_{r+2}-f_{r+2})\cap(A_{r+2}'-f_{r+2}')}
{\#((F_{r,r+2}-f_{r+2})\cap(F_{r,r+2}'-f_{r+2}')}.
$$
Slightly modifying the argument in the proof of Theorem~5.2 we deduce \thetag{6-2}.
If either the second or the third condition from \thetag{6-3} is satisfied for infinitely many $n$ then one can easily deduce  from \thetag{6-2} and \thetag{6-3} that one of the coordinate projection of $\Cal G(\lambda)$ is of 0 measure.
If the second and the third condition  are satisfied for only finitely many $n$ then the first condition in \thetag{6-3} is satisfied for all but finitely many $n$.
Hence it is satisfied for infinitely many elements of $J$.
However, if $\widetilde c_n=\widetilde c_n'$ for $n\in J$ then $\widetilde c_n=\widetilde c_n'=0$.
Therefore the two marginals of $\Cal G(\lambda)$ are of 0 measure.
\qed
\enddemo

We also note that  one can construct\footnote{Via a slight modification of Example~5.4.}, for each $\omega$, a non-ergodic Radon 2-fold self-joining of $T^\omega_1$ whose ergodic components are not joinings because their coordinate projections are singular to $\mu^\omega$.

Given $\omega\in\{0,1\}^\Bbb N$, we define $\omega^*\in\{0,1\}^\Bbb N$ by setting $\omega^*(i)=1-\omega(i)$.
It is easy to see that $T_1^{\omega^*}$ is conjugate to the inverse to $T_1^\omega$.

\proclaim{Corollary 6.2}  For each $\omega$, the transformation $T_1^\omega$ is not conjugate to its inverse.
Moreover, $T_1^\omega$ and its inverse are  Radon disjoint.
\endproclaim

\head 7. Further generalizations and some open problems
\endhead

It is easy to see that the above argument works almost verbally for a more general class of rank-one transformations.
To specify this class we first give an alternative (but equivalent) description of the infinite Chacon transformation.
For that we will use the classical language of cutting-and-stacking construction \cite{Fr}.
The initial 0-th tower consists of a single interval $[0,1)$.
We now describe the inductive procedure of  passing from the $n$-th tower consisting of $h_n$ levels of width $\frac 1{3^n}$ to the $(n+1)$-th tower.
For that we cut the $n$-th tower into 3 subtowers (called copies) of equal width.
Then we place the second copy  over the first one, add an additional level (called spacer) over the second copy and put the third copy over this spacer.
Next we put $[1.5h_n]$ spacers over the top of the third copy and
$[1.5h_n]+1$ spacers under the bottom of the first copy.
Here $[.]$ stands for the integer part.
We thus obtain the $(n+1)$-tower consisting of $h_{n+1}=6h_n+1$ levels of width $\frac 1{3^{n+1}}$.
The  transformation $T$ moves each (except for the highest  one) level of the tower one level up.
The transformation is not defined on the highest level of the tower.
However in the limit we obtain a well defined transformation on $\Bbb R$ (which is the union of all levels of all towers) endowed with Lebesgue measure.
It is measure theoretically isomorphic to
 the infinite Chacon transformation described above via the $(C,F)$-construction in Section 4.

We now can obtain a family of transformations using almost the same cutting-and-stacking algorithm but adding ``more'' spacers.
I mean the following.
Let $(\alpha_n)_{n\in\Bbb N}$ and $(\beta_n)_{n\in\Bbb N}$ be two arbitrary sequences of nonnegative integers.
When constructing  the $(n+1)$-th tower we first do verbally what we did in the above construction of  the infinite Chacon transformation and after that we put $\alpha_n$ additional spacers on the top and $\beta_n$ spacers under the bottom of the tower.
In the limit of the inductive construction we obtain a
 transformation that  has  Radon MSJ.
The proof is almost the same as for the infinite Chacon transformation.

We conclude the paper with a list of open problems.
\roster
  \item"---" Are there two ergodic Radon dynamical systems which are Radon disjoint but 
  which admit a non-ergodic Radon joining?
  \item"---" Does the infinite Chacon transformation have  invariant sub-$\sigma$-algebras admitting equivalent $\sigma$-finite invariant measures?
   \item"---" More generally, given a $(C,F)$-transformation $T$ with an infinite invariant Radon measure $\mu$ such that all ergodic  $T\times T$-invariant Radon measures are either quasi-graphs or $\mu\times\mu$, does it have factors admitting equivalent $\sigma$-finite invariant measures?
\endroster

\enddocument